  \documentclass[12pt]{amsart}
  \usepackage{latexsym} 
  \usepackage[all]{xy}
  \usepackage{amsfonts} 
  \usepackage{amsthm} 
  \usepackage{amsmath} 
  \usepackage{amssymb}
  \usepackage{pifont}

  \def\sw#1{{\sb{(#1)}}}

  \def\suc#1{{\sp{(#1)}}} 
   
  \def\tens{\mathop{\otimes}}

  \def\<{{\langle}} 
  \def\>{{\rangle}}

  \def\eps{\varepsilon}

  \def\note#1{{}} 
 \def\can{{\rm \textsf{can}}} 

\def\cocan{{\overline{\rm \textsf{can}}}}

  \def\cocan{\overline{\can}} 
  \def\note#1{}

  \def\cA{{\mathcal A}}

  \def\beq{\begin{equation}} 
  \def\eeq{\end{equation}}

  \def\ot{{\otimes}}

  \def\Aro{{{}^A\!\varrho}}
  \def\roA{{\varrho^A}}

 \def\ta{\tilde{a}}

 \def\re{\psi}

  \newcounter{zlist}

  \newcounter{blist} 
  \newenvironment{blist}{\begin{list}{(\alph{blist})}{ 
  \usecounter{blist}\leftmargin2.5em\labelwidth2em\labelsep0.5em 
  \topsep0.6ex 
  \parsep0.3ex plus0.2ex minus0.1ex}}{\end{list}} 

  \newcounter{rlist}

\def\stac#1{\raise-.2cm\hbox{$\stackrel{\displaystyle\otimes}{\scriptscriptstyle{#1}}$}}

\def\cten#1{\raise-.2cm\hbox{$\stackrel{\displaystyle\widehat{\otimes}}
{\scriptscriptstyle{#1}}$}}

  \def\Label#1{\label{#1}\ifmmode\llap{[#1] }\else 
  \marginpar{\smash{\hbox{\tiny [#1]}}}\fi} 
  \def\Label{\label}

  \newtheorem{proposition}{Proposition}

  \newtheorem{theorem}[proposition]{Theorem} 

  \theoremstyle{definition} 
  \newtheorem{definition}[proposition]{Definition}

  \theoremstyle{remark}

  \newcounter{c} 
   
  \newcommand{\etyk}[1]{\vspace{-7.4mm}$$\begin{equation}\Label{#1} 
  \addtocounter{c}{1}} 
  \renewcommand{\]}{\ifnum \value{c}=1 $$\else \end{equation}\fi} 
\begin{document} 

 \title{An explicit formula for a strong connection} 
 \author{E.J.\ Beggs \& Tomasz Brzezi\'nski}
 \address{ Department of Mathematics, Swansea University, 
  Singleton Park, \newline\indent  Swansea SA2 8PP, U.K.} 
  \email{E.J.Beggs@swansea.ac.uk}
  \email{T.Brzezinski@swansea.ac.uk}   
    \date{June 2006}
  \begin{abstract} 
 An explicit formula for a strong connection form in a principal extension by a coseparable coalgebra is given.
   \end{abstract} 
  \maketitle 

\noindent {\bf 1.} In the studies of geometry of non-commutative principal bundles or coalgebra-Galois extensions (cf.\ \cite{BrzMaj:coa}) an important role is played by the notion of a {\em strong connection} (for the universal differential structure) first introduced in the context of Hopf-Galois extensions in \cite{Haj:str}. The existence of a strong connection guarantees that a bundle associated to a coalgebra-Galois extension is a (finitely generated and) projective module, hence it is (a module of sections on) a vector bundle in the sense of non-commutative geometry (cf.\ \cite{Con:ncg}). Furthermore, a strong connection form gives rise to a {\em Chern-Galois character} \cite{BrzHaj:che}, a mapping from the Grothendieck group of isomorphism classes of finite dimensional corepresentations of the structure coalgebra to the cyclic homology of the base algebra (see \cite{BohBrz:str} for the most general, relative formulation).

The existence of a strong connection in a Hopf-Galois extension is assured by the classical Schneider Theorem I \cite{Sch:pri}. This states that a free coaction of a Hopf algebra with bijective antipode on its injective comodule algebra determines a Hopf-Galois extension with a strong connection. This theorem has been extended to coalgebra (entwined) extensions with a {\em coseparable colagebra} \cite[Theorem~4.6]{Brz:gal} \cite[Theorem~5.9]{SchSch:gen}. In all these cases the proof of existence is not a constructive proof: the existence follows by general arguments, but no explicit form of connection is given. On the other hand, the knowledge of this form is needed for construction and calculation of Levi-Civit\'a connections and projectors for associated bundles, and the  Chern-Galois characters. Recently, several examples of strong connections have been constructed (cf.\ \cite{BrzMaj:geo}, \cite{HajMat:loc},  \cite{Lanvan:pri}, \cite{LanPag:Hop}) or their form conjectured \cite{BonCic:bij}, but no general procedure has been established. The aim of this note is to give a direct proof of a Schneider type theorem for coalgebra extensions in which the connection is explicitly given.

We work over a field $k$, unadorned tensor product is over $k$. For a vector space $V$, the identity map is denoted by the same symbol $V$. All algebras are associative and unital. In an algebra $A$, $1$ denotes the unit both as an element and as a $k$-linear map $k\to A$ and $\mu:A\ot A\to A$ denotes the product. In a coalgebra $C$, the coproduct is denoted by $\Delta$ and counit by $\eps$. We denote coactions of $C$ on a vector space $A$ by $\roA$ (the right coaction) and $\Aro$ (the left coaction). The following Sweedler's notation is used: $\Delta(c) = c\sw 1\ot c\sw 2$, $\roA(a) = a\sw 0\ot a\sw 1$, $\Aro(a) = a\sw{-1}\ot a\sw 0$ (summation suppressed).  \bigskip

\noindent {\bf 2.} In this section a strong connection form is explicitly constructed from a cointegral in a coseparable coalgebra. First we recall the definition of a non-commutative object which captures most of the geometric information carried by (locally trivial) principal bundles.
\begin{definition}[cf.\ \cite{BrzHaj:che}, Definition~2.1]
\label{def.principal}
Let $C$ be a coalgebra and $A$ an algebra and a right $C$-comodule via
$\varrho^A:A\to A\ot C$. Let
$$
B=A^{co C}:=\{b\in A~|~\varrho^A(ba)=b\varrho^A(a),\ \forall a\in A\},
$$
 denote the subalgebra of $C$-coinvariants of $A$. 
The inclusion of algebras $B\subseteq A$ is called a
 {\em principal $C$-extension}
 if
\begin{blist}
\item
$
\can: A\ot_BA{\to} A\ot C,\;a\ot a'\mapsto a\varrho^A(a')
$
is bijective (the Galois  condition);
\item 
$A$ is $C$-equivariantly projective as a left $B$-module, i.e.\ there exists a left $B$-module, right $C$-comodule section of the product $B\ot A\to A$;
\item  $\psi:C\ot A{\to} A\ot C$, $c\ot a\mapsto
\can(\can^{-1}(1\ot c)a)$ is bijective;
\item 
there is a group-like element $e\in C$ such that $\varrho^A(a)=\psi(e\ot a)$,
for all $a\in
A$.
\end{blist}
\end{definition}

By  \cite[Theorem~3.5]{BrzHaj:coa}, the map $\psi$ defined in Definition~\ref{def.principal}(c) is an example of a {\em right-right entwining map}, i.e.\ it is a map
$\re : C \otimes A \to A \otimes C$, 
which, 
    satisfies
    the following relations:
    $$
    \psi\circ({C}\tens \mu) = (\mu\tens {C})\circ
       (A\tens\psi)\circ(\psi\tens
 {A}), \qquad 
\psi\circ ({C}\tens 1) = 1\tens {C},
 $$
 $$
({A}\tens\Delta)\circ\psi = (\psi\tens
      C)\circ({C}\tens\psi)\circ(\Delta\tens {A}), \qquad
({A}\tens \eps)\circ\psi =
\eps\tens A.
$$
Consequently, the inverse of $\psi$ is a {\em left-left entwining map}, i.e.\ 
    the following relations
  $$
    \psi^{-1}\circ(\mu\tens C) = (C\tens \mu)\circ
      (\psi^{-1}\tens
 {A})\circ  (A\tens\psi^{-1}), \qquad 
{C}\tens 1 = \psi^{-1}\circ (1\tens {C}),
 $$
 \begin{equation}\label{le2}
(\Delta\tens)\circ\psi^{-1} = ({C}\tens\psi^{-1})\circ(\psi^{-1}\tens
      C)\circ(A\tens \Delta), \qquad
{A}\tens \eps =
(\eps\tens A)\circ\psi^{-1}.
\end{equation}
are satisfied. 

Furthermore, $A$ is a {\em right entwined module}, i.e.\ the map $\psi$ makes the $C$-coaction $\roA$ compatible with the product in the sense that, for all $a,\ta\in A$,
\begin{equation}
\roA(a\ta) = a\sw 0\psi (a\sw 1\ot \ta).
\label{rem}
\end{equation}
Since $\psi$ is bijective, $A$ is also a left $C$-comodule with the coaction 
\begin{equation}\label{gen.left.coa}
\forall a\in A, \quad \Aro(a)=\psi^{-1}(a\roA(1)).
\end{equation}
In view of condition (d) in Definition~\ref{def.principal}, in the case of a principal extension this left coaction comes out explicitly as
$\roA (a) = \psi^{-1}(e\ot a)$.
With this coaction $A$ is a {\em left entwined module}, i.e., for all $a, \ta\in A$, $\Aro(a\ta) = \psi^{-1}(a\ot \ta\sw{-1})\ta\sw 0$. Note that if $C$ is a Hopf algebra with a bijective antipode $S$ and $A$ is a right $C$-comodule algebra, then the map $\psi$ in Definition~\ref{def.principal}(c) and its inverse come out as
$$
\psi (c\ot a) = a\sw 0 \ot ca\sw 1, \qquad \psi^{-1}(a\ot c) = cS^{-1}a\sw 1\ot a\sw 0.
$$
Hence, by setting $e=1_C$ we obtain that $\roA(a) = \psi( 1_C\ot a)$ and $\Aro(a) = S^{-1}a\sw 1\ot a\sw 0$. In general, we make the following 
\begin{definition}[cf.\  Definition~2.2 in \cite{BohBrz:str}]\label{def.ee}
Let $C$ be a coalgebra and let $A$ be an algebra and a right $C$-comodule. An inclusion of algebras $B\subseteq A$ is called an {\em entwined $C$-extension} if $B$ is a subalgebra of coinvariants $B=A^{coC}$ and there exists a bijective right-right entwining map $\psi: C\ot A\to A\ot C$ such that compatibility condition \eqref{rem} is satisfied.
\end{definition}
Note that the compatibility condition \eqref{rem} imply that the right coaction in an entwined extension is given by
\begin{equation}\label{gen.right.coa}
\forall a\in A, \quad \roA(a)= 1\sw 0\psi(1\sw 1\ot a).
\end{equation}

In particular, conditions (a) and (c) in Definition~\ref{def.principal} imply that a principal extension is an example of an entwined extension. The condition (c) in Definition~\ref{def.principal} is equivalent to the existence of a {\em strong connection form}.
\begin{definition} \label{str.con.form}
Let $B\subseteq A$ be an entwined $C$-extension and let $\cocan$ be  the {\em lifted canonical map},
$$
\cocan: A\ot A{\to} A\ot C,\;a\ot a'\mapsto a\varrho^A(a').
$$
A $k$-linear map $\ell: C\to A\ot A$ satisfying the following properties:
\begin{blist}
\item  $\cocan\circ\ell = 1_A\otimes C$;
\item $(\ell\otimes C)\circ\Delta= (A\otimes \roA)\circ\ell$;
\item $(C\otimes\ell)\circ\Delta = (\Aro\otimes A)\circ\ell$,
\end{blist}
is called a {\em strong connection form} or a {\em strong connection lifting}.
 Here $\Aro$ is the induced left coaction as in \eqref{gen.left.coa}.
\end{definition}
 Existence of a strong connection form $\ell$ in an entwined $C$-extension $B\subseteq A$  implies that it is a Galois extension, i.e.\ the canonical map $\can$ is bijective, and that $A$ is $C$-equivariantly projective as a left $B$-module (cf. \cite[Theorem~3.7, Corollary~3.8]{BohBrz:str} for a detailed proof in the most general case). Explicitly, the splitting $s: A\to B\ot A$ of the product is given by $s(a) = a\sw 0\ell(a\sw 1)$. If, in addition, there is a group-like element $e\in C$ such that $\roA(1) =1\ot e$, then an entwined extension with a strong connection form is a principal extension. In this case, $\ell$ can always be chosen in such a way that $\ell(e) = 1\ot 1$. The existence of a group-like element $e$ is needed in order to have a bijective correspondence between strong connection forms $\ell$ and strong covariant differentials (and also to make the universal differential calculus on $A$ a {\em $C$-covariant} calculus, cf.\ discussion in \cite[Sections~4,5]{BrzMaj:geo}).

Recall that a coalgebra $C$ is said to be {\em coseparable} if the coproduct has a retraction in the category of $C$-bicomodules, equivalently, if there exists a $k$-linear map $\delta: C\ot C\to k$ with the following properties, for all $c,c'\in C$,
\begin{equation}
\delta(c\sw 1\ot c\sw 2) = \eps(c), \qquad c\sw 1\delta(c\sw 2\ot c') = \delta(c\ot c'\sw 1)c'\sw 2.
\label{coint}
\end{equation}
Such a map $\delta$ is called a {\em cointegral}. For example, if $C$ is a coalgebra spanned by a set of group-like elements $x_i$, then $C$ is a coseparable with a cointegral
$
\delta(x_i\ot x_j) = \delta_{ij}.
$
If $C$ is a Hopf algebra, then it is coseparable if and only if there exists a normalised left (or right) integral on $C$, i.e.\ a linear map $\lambda : C\to k$ such that, for all $c\in C$,
$$
c\sw 1\lambda(c\sw 2)= \lambda(c), \qquad \lambda(1) =1.
 $$
 The corresponding cointegral is
 $
 \delta (c\ot c')  = \lambda \left(c S(c')\right).
 $
Conversely, given a cointegral $\delta$ on a Hopf algebra $C$, the left integral is obtained as $\lambda( c )= \delta(c\ot 1)$. Since, by the Woronowicz theorem \cite[Theorem~4.2]{Wor:com} every compact quantum group has an integral (or a {\em Haar measure}) most of the 
 coalgebras which are of interest in non-commutative differential geometry are coseparable.

The main result of this note is contained in the following theorem, which gives the explicit form of a strong connection.
\begin{theorem}
\label{thm.main}
Let $B\subseteq A$ be an entwined $C$-extension. Assume that $C$ is a coseparable coalgebra with a cointegral $\delta: C\ot C\to k$ and that the (lifted) canonical map
$$
\cocan: A\ot A{\to} A\ot C,\;a\ot a'\mapsto a\varrho^A(a')
$$
is surjective. Write  $\sigma: C\to A\ot A$ for a $k$-linear map such that $\cocan\circ\sigma = 1\ot C$ and define maps $\gamma : C\ot A\to A$ and $\alpha : A\ot C\to A$ by
$$
\gamma = (\delta\ot A)\circ (C\ot \Aro), \qquad \alpha = (A\ot \delta)\circ (\roA\ot C).
$$
Then 
\begin{equation}
\label{ell}
\ell = (\gamma\ot\alpha)\circ (C\ot\sigma\ot C)\circ (\Delta\ot C)\circ \Delta, 
\end{equation}
is a strong connection form. 

Furthermore, if $\roA(1) = 1\ot e$, then $A$ is a principal $C$-extension.
\end{theorem}
\begin{proof}
Using the definition of a cointegral, one easily checks that the map $\gamma$ is left $C$-colinear, where $C\ot A$ as understood as a left $C$-comodule via $\Delta\ot A$, and $\alpha$ is right $C$-colinear, where $A\ot C$ is a right $C$-comodule via $A\ot \Delta$. By the colinearity of $\gamma$ and $\alpha$ the map $\ell$ is $C$-bicolinear. 

To prove that $\ell$ is a section of the map $\cocan$ we start with the following simple calculation, for all $a,\ta\in A$,
$$
\psi^{-1}(a\ta\sw 0\ot \ta\sw 1) = 
 \psi^{-1}\left(a1\sw 0\psi(1\sw 1\ot\ta)\right)
= \psi^{-1}(a1\sw 0\ot 1\sw 1)\ta = a\sw{-1}\ot a\sw 0\ta.
$$
Here the first and last equalities follow from the definitions of the right and left $C$-coactions on $A$ (cf.\ \eqref{gen.right.coa},  \eqref{gen.left.coa}), and the second equality follows by \eqref{le2} and by the fact that $\psi^{-1}$ is the inverse of $\psi$. Thus we obtain the equality
\begin{equation}
\psi^{-1}(a\ta\sw 0\ot \ta\sw 1)\ot \ta\sw 2 = a\sw{-1}\ot a\sw 0\ta\sw 0\ot \ta\sw 1. \label{key}
\end{equation}
For any $c\in C$,  write explicitly
$
c\suc 1 \ot c\suc 2 := \sigma(c),
$
so that 
$
c\suc 1c\suc 2\sw 0\ot c\suc 2\sw 1 = 1\ot c.
$
This leads to the equality
$$
c\sw 1\ot c\sw 2\suc 1c\sw 2\suc 2\sw 0  \ot c\sw 2\suc 2\sw 1\ot c\sw 3 = c\sw 1\ot 1\ot c\sw 2 \ot c\sw 3.
$$
Apply $(C\ot \psi^{-1}\ot C\ot \Delta)\circ (C\ot A\ot\Delta\ot C)$ and then use \eqref{key} on the left hand side and \eqref{le2} on the right hand side to obtain
$$
c\sw 1\ot c\sw 2\suc 1\sw{-1} \ot c\sw 2\suc 1\sw{0}c\sw 2\suc 2\sw 0 
 \ot \, c\sw 2\suc 2\sw 1\ot c\sw 3\ot c\sw 4
 = c\sw 1\ot c\sw 2\ot 1\ot c\sw 3 \ot c\sw 4\ot c\sw 5.
$$
Now apply $\delta\ot A\ot\delta\ot C$ and use the definitions of maps $\gamma$ and $\alpha$ in terms of $\delta$ on the left hand side, and the properties of the cointegral \eqref{coint} on the right, to conclude that
$$
\gamma(c\sw 1 \ot c\sw 2\suc 1)\alpha(c\sw 2\suc 2\ot c\sw 3)\ot c\sw 4 = 1\ot c.
$$
By the right $C$-colinearity of $\alpha$ this implies that $\cocan \circ \ell = 1\ot C$ as required. In view of the discussion that follows Definition~\ref{str.con.form}, the second assertion is obvious.
\end{proof}

Note that if $\roA(1) = 1\ot e$, then  $\cocan (1\ot 1) = 1\ot e$, hence the linear map $\sigma$ can always be normalised so that $\sigma(e) = 1\ot 1$ by making the linear change
$$
\sigma \mapsto \sigma +1\ot 1\eps - \sigma(e)\eps.
$$
The strong connection form obtained with such normalised $\sigma$ is also normalised, i.e.\ $\ell(e) = 1\ot 1$. 
Furthermore, if $\sigma$ is right (resp.\ left) $C$-colinear, then the formula \eqref{ell} reduces to
$$
\ell = (\gamma \ot A)\circ (C\ot\sigma)\circ\Delta, \qquad \mbox{(resp.\ 
$\ell = (A\ot \alpha)\circ (\sigma\ot C)\circ\Delta$)}.
$$
Thus, if $\sigma$ is $C$-bicolinear, then $\ell =\sigma$. This gives an effective way of testing whether the map $\sigma$ is bicolinear. \bigskip

\noindent {\bf 3.} The main usefulness of formula \eqref{ell} lies in the fact that usually the map $\sigma$ in Theorem~\ref{thm.main} is already obtained as the first step of checking whether a given extension is a Galois extension. Furthermore, in geometrically most interesting cases, the coalgebra $C$ is a Hopf algebra (although the coaction is not necessarily an algebra map), for which the explicit form of the left integral (or Haar measure) is known. As an illustration one can consider an entwined extension 
$\cA(\Sigma_q^4)\subseteq \cA(S_q^7)$ by a coalgebra $\cA(SU_q(2))$ constructed in \cite{BonCic:ins}. Here $\cA(\Sigma_q^4)$ is the algebra of functions on a quantum four-sphere, $ \cA(S_q^7)$ is the algebra of functions on a quantum seven-sphere (obtained as a quotient of the quantum group $U(4)$), while $\cA(SU_q(2))$ is the algebra of functions on the quantum group $SU(2)$. The $k$-linear splitting $\sigma$ required in Theorem~\ref{thm.main} is constructed explicitly in \cite[Equations~(11)]{BonCic:bij} (as a part of a complete proof that this is a coalgebra-Galois extension), while the formula for the left integral on $\cA(SU_q(2))$ is derived  in \cite[Appendix A1]{Wor:com}. Once these two are combined with each other, finding a strong connection $\ell$ \eqref{ell} is a matter of straightforward albeit tedious calculations. In view of the discussion at the end of  {\textsection 2}, this provides one also with a means of checking if  \cite[Equations~(11)]{BonCic:bij} indeed define a strong connection form as conjectured at the end of \cite{BonCic:bij}. 

An example of quantum principal bundles is provided by {\em quantum homogeneous spaces}. In this case $A$ is a Hopf algebra and a quantum homogeneous space is defined as a subalgebra $B\subseteq A$ such that $\Delta(B) \subseteq A\ot B$. The coalgebra $C$ is defined as a quotient $C = A/B^+A$, where $B^+ = B\cap\ker\eps$, while  the surjection $\pi: A\to C$ induces left and right coactions of $C$ on $A$ via $(\pi\ot A)\circ\Delta$ and $(A\ot\pi)\circ\Delta$. If the antipode in $A$ is bijective, this produces an entwined extension $A^{coC}\subseteq A$ (and $B=A^{coC}$ for example if there is a strong connection). By \cite[Proposition~4.4]{BrzMaj:geo}, left invariant strong connection forms are in bijective correspondence with $C$-bicolinear maps $\iota: C\to A$ such that $\pi\circ\iota = C$. In view of Theorem~\ref{thm.main} (or directly using \eqref{coint}), if $C$ is a coseparable coalgebra with a cointegral $\delta$, then any $k$-linear section $i$ of $\pi$ gives rise to such a $C$-bicolinear section $\iota$ by the formula, for all $c\in C$,
\begin{equation}
\iota(c) = \delta\left( c\sw 1 \ot \pi(i(c\sw 2)\sw 1)\right)i(c\sw 2)\sw 2\delta\left( \pi(i(c\sw 2))\sw 3\ot c\sw 3\right).
\label{iota}
\end{equation}
For example, \cite[Proposition~6.1]{BrzMaj:geo} gives an explicit formula for such a $k$-linear section of $\pi$ in the case of a quantum Hopf fibration over a general quantum two-sphere and hence a $C$-bicolinear map such as in \cite[Proposition~6.3]{BrzMaj:geo} can be obtained by the above averaging procedure \eqref{iota}.

Finally, the explicit formula \eqref{ell} can be used to calculate Chern-Galois characters for principal extensions with coseparable coalgebras. The components of the Chern-Galois character are
defined in \cite[Corollary~3.2]{BrzHaj:che} in terms of the strong connection form $\ell$. However, settling the question whether the explicit knowledge of $\ell$ in terms of a cointegral provides one with an effective method of gaining information about a principal extension seems to require a case by case analysis. This is hoped to be attempted elsewhere.

\end{document}